\documentclass[12pt]{amsart}

\usepackage{amsmath,amssymb,latexsym,enumerate}
\usepackage{a4wide}
\usepackage{color,graphicx}
\usepackage{ifthen}
\usepackage{verbatim}
\usepackage{fancyhdr}
\usepackage{url}
\usepackage{bbm}
\usepackage{graphicx}
\usepackage[sans]{dsfont}
\usepackage{mathrsfs}
\usepackage{MnSymbol}

\newcommand{\Id}{\operatorname{Id}}
\newcommand{\Op}{\operatorname{Op}}

\newcommand{\supp}{\operatorname{supp}}

\newcommand{\CC}{\mathscr{C}}
\newcommand{\N}{\mathbb N}
\newcommand{\R}{\mathbb R}

\def\<{\langle}
\def\>{\rangle}
\newcommand{\bp}{{\it Proof. }}
\newcommand{\ep}{\hfill $\square$\\}

\newcommand{\be}{\begin{equation}}
\newcommand{\ee}{\end{equation}}
\newcommand{\bes}{\begin{equation*}}
\newcommand{\ees}{\end{equation*}}

\renewcommand{\dh}{d_{\varphi,h}}

\newcommand{\dhs}{d_{\varphi,h}^{G,*}}

\numberwithin{equation}{section}
\numberwithin{figure}{section}

\newtheorem{theorem}{Theorem}[section]

\newtheorem{lemma}[theorem]{Lemma}
\newtheorem{defin}[theorem]{Definition}
\newtheorem{proposition}[theorem]{Proposition}
\newtheorem{remark}[theorem]{Remark}

\def\ccc{{\mathcal C}}
 
 \def\lll{{\mathcal L}}
\def\mmm{{\mathcal M}} \def\ooo{{\mathcal O}}

\def\vvv{{\mathcal V}}

\begin{document}
\title{Around  supersymmetry for semiclassical second order differential operators}

\author[L. Michel]{Laurent~Michel}
\address{Laboratoire J.-A. Dieudonn\'e\\ Universit\'e de Nice}
\email{lmichel@\allowbreak unice.\allowbreak fr}

\begin{abstract}
Let $P(h),h\in]0,1]$ be a semiclassical scalar differential operator of order $2$. The existence of  a supersymmetric structure given by a matrix $G(x;h)$ 
was  exhibited in \cite{HeHiSj13} under rather general assumptions. In this note we 
 give a sufficient condition the coefficients of $P(h)$ so that the matrix $G(x;h)$ enjoys some nice estimates with respect to 
 the semiclassical parameter.

\end{abstract}
\maketitle

\section{Introduction}
In many problems arising in physics one is interested in computing accurately the spectrum of some differential operators depending on
a small parameter (that we shall denote by $h$ throughout). In numerous situations, proving sharp results 
can be done by using specific structure of the operator. For instance, in the setting of Schr\"odinger operators, 
geometric assumptions on the potential leads to sharp computation of the splitting between eigenvalues \cite{HeSj84_01}.
More recently, the computation of the low lying eigenvalues of 
the semiclassical Witten Laplacian was performed by using the specific structure of the operator \cite{HeSj85_01}, \cite{HeKlNi04_01}. 
In these papers, the fact that the Witten Laplacian enjoys a supersymmetric structure (that is can be written as a twisted Hodge Laplacian) is fundamental and 
it doesn't seem possible to obtain the sharp results without using this property.
Similarly, existence of supersymmetric structure was used in \cite{BoHeMi15} to compute the spectrum of some semiclassical Markov operator and hence the rate of convergence to equilibrium of the associated random walk.

In a nonselfadjoint setting numerous results in the same spirit were obtained by H\'erau, Hitrik, Sj\"ostrand \cite{HeHiSj08-2,HeHiSj08-1,HeHiSj11_01}. 
In all these papers, the authors are lead to compute a spectral gap for non self-adjoint operators. 
Their approach is based on the fact that the underlying operator is supersymmetric for a convenient bilinear product and then some tools developed for the study
of Witten Laplacian can be used (of course, one major additional difficulty comes from the fact that they are in a nonselfadjoint situation).

In most situations mentioned above, the supersymmetric structure of the operator is known in advance. 
Nevertheless, it can occurs that  the supersymmetric structure is hidden and has to be exhibited. 
This was for instance the case in \cite{BoHeMi15} where the authors give a sufficient condition for selfadjoint 
pseudodifferential operators to be supersymmetric. Roughly speaking, the main  assumption made in \cite{BoHeMi15} 
is that the Weyl symbol of the operator is an even function of the $\xi$ variable. The first motivation of the present  work, was hence to investigate 
what happends when this assumption fails to be true. As we shall see later, 
the pseudodifferential situation is quite intricate and we shall restict our attention to the case of second order scalar differential operators.
This issue was already addressed in \cite{HeHiSj13}, where the authors consider 
semiclassical and nonsemiclassical operators $P$. 
In both situation the author exhibit supersymmetric structure (under the assumption that the kernels of $P$ and $P^*$ contain some specific element) but they emphasize the fact that in the semiclassical situation, 
the factorization of the operator is done without control with respect to the semiclassical parameter. 
The goal of this  paper is to give a sufficient condition in order to have a control with respect to $h$ in the 
factorization and also to discuss the optimality of this condition. We decided to use here the formalism of \cite{HeHiSj13} that we recall in the next paragraph.

Let $X$ be either $\R^n$ either a $n$-dimensional smooth connected compact manifold without boundary equipped with a smooth volume density $\omega (dx)$, and 
let $P=P(h)$, $h\in]0,1]$ denote a second order scalar semiclassical differential operator on $X$ with real smooth coefficient. 
For $0\leq k<n$, let $\Omega^k(X)=\ccc^\infty(X,\Lambda^kT^*X)$ and denote $d:\Omega^k(X)\rightarrow \Omega^{k+1}(X)$  the exterior derivative. For any $x\in X$, we recall the natural pairing $\<.,.\>_{\Lambda,\Lambda^*}$ on 
$\Lambda^kTX\times\Lambda^kT^*X$ given by $\<u,v^*\>_{\Lambda,\Lambda^*}=\det ((v^*_i(u_j)_{i,j}))$ for any $v^*=v_1^*\wedge\ldots\wedge v_k^*$ and $u=u_1\wedge\ldots\wedge u_k$. It gives rise to a natural pairing on $\ccc^\infty(X,\Lambda^kT^*X)\times\ccc^\infty(X,\Lambda^kTX)$ by integrating the preceding formula against the volume form. Then, we let 
$\delta: \ccc^\infty(X,\Lambda^kTX)\rightarrow \ccc^\infty(X,\Lambda^{k-1}TX)$ be the adjoint of $d$ for this pairing.

Suppose that $G(x): T_x^*X\rightarrow T_xX$ is a linear mapping depending smoothly on $x\in X$. 
Then $\Lambda^kG$ maps $\Lambda^kT_x^*X$ into $\Lambda^kT_xX$ (by convention $\Lambda^0G$ is the identity on $\R$) 
and we can define a bilinear product on $\ccc_c^\infty(X,\Lambda^kT^*X)$ by the formula
\be\label{eq:def_barcketG}
\<u,v\>_G=\int_X\<G(x)u(x),v(x)\>_{\Lambda,\Lambda^*}\omega(dx).
\ee
where for short, we write $G(x)$ instead of $\Lambda^kG(x)$.
When $G(x)$ is invertible for any $x\in X$ we can define 
$d^{G,*}=(G^t)^{-1}\delta G^t$ and one checks easily that $d^{G,*}$ is the formal adjoint 
of $d$ with respect to $G$ 
\be
\<d u,v\>_G=\<u,d^{G,*}v\>_G,\;\forall u,v\in \ccc_c^\infty(X,\Lambda^kT^*X)
\ee
Notice that on $1$-forms $d^{G,*}: \ccc^\infty(X,\Lambda^1T^*X)\rightarrow \ccc^\infty(X,\R)$ is given by $d^{G,*}=\delta \circ G^t$ 
which makes sense even if $G$ is not invertible.

In the case where $X$ is a compact Riemaniann manifold, we can identify $T_x^*X $ and $T_xX$ by mean of the metric $g$, so that  $G$ can be considered 
as an operator acting on $T_x^*X$. When $X=\R^n$ is equipped with the euclidean metric, 
then $G$ will be identified with its matrix in the basis of canonical $1$-forms.

Given $\varphi\in \ccc^\infty(X,\R)$, the associated Witten complex is defined by the semiclassical weighted de Rham differentiation
$$\dh=e^{-\varphi/h}\circ h d\circ e^{\varphi/h}=hd+d\varphi^\wedge$$
and its formal adjoint with respect to the bilinear form \eqref{eq:def_barcketG} is
$$\dhs=e^{\varphi/h}\circ h d^{G,*}\circ e^{-\varphi/h}=(G^t)^{-1}\circ(h\delta+d\varphi^\righthalfcup)\circ G^t.$$
Let us now recall the definition of supersymmetric structure used in \cite{HeHiSj13}
\begin{defin}
Let $P=P(x,hD_x;h)$ be a second order scalar real semiclassical differential operator on $X$. We say that $P$ has a supersymmetric structure if there exists a linear $h$-dependent map 
$G(x;h):T_x^*X\rightarrow T_xX$, 
smooth with respect to  $x\in X$ and functions $\varphi,\psi\in \ccc^\infty(X,\R)$ such that 
$$P=d_{\psi,h}^{G,*}d_{\varphi,h}$$
for all $h\in]0,h_0],h_0>0$.
\end{defin}
Here we decided to consider phase functions $\varphi, \psi$ which are independent of $h$ in order to simplify.
As noticed in \cite{HeHiSj13}, no control of $G(x;h)$ with respect to $h$ is required in this definition. 
In order to get some bounds on  $G(x;h)$, we first need to handle a metric on $X$. If $X=\R^n$ we consider $g$ the Euclidean metric and if 
$X$ is a compact manifold we take $g$ to be any Riemaniann metric. From this metric we get normed vector space structure on $T_xX$ and $T^*_xX$. 
In the case where $X=\R^n$ we need to control the function at infinity. Given an order function $a$ (in the sense of Def. 7.4 in \cite{DiSj99_01}) 
we say that a 
function $f\in \ccc^\infty(X,\R)$ belongs to $S(a)$ if 
\be\label{eq:def-sym}
\forall \alpha\in \N^n,\exists C_\alpha>0,\,\forall x\in X,\;|\partial^\alpha f(x)|\leq C_\alpha a(x)
\ee
Throughout, given $m\in\R$, we will often use the order function $\rho_m$ defined 
by $\rho_m(x)=(1+|x|^2)^{\frac  m 2}$ when $X=\R^n$ and by $\rho_m=1$ if $X$ is compact.
We shall denote $S_m=S(\rho_m)$.
 We introduce the following
\begin{defin}
Let $P=P(x,hD_x;h)$ be a second order scalar real semiclassical differential operator on $X$. 
We say that $P$ has a temperate supersymmetric structure if it has a supersymmetric structure (in the sense of the above definition) and if 
the map $G(x;h):T_x^*X\rightarrow T_xX$ satisfies the following:
 there exist $m\in\R$ and some constant $C_{\nu}>0$ such that 
\be
\label{eq:def-temp}
 \|\partial_x^\nu G(x,h)\Vert_{T_x^*X\rightarrow T_xX}\leq C_{\nu}\rho_m(x),\;\forall x\in X
\ee
for all $h\in]0,h_0]$. 
\end{defin}
 Throughout the paper we shall call ``supersymmetric structure'' any operator $G(x;h)$ as above.
 We shall say that $G$ is temperate if it satisfies \eqref{eq:def-temp} for some $m\in\R$.
Observe that the preceding definition doesn't depend on the choice of the metric $g$ since if $g_1$ and $g_2$ are 
two metric on a compact manifold $X$, the corresponding norms on tangent and cotangent spaces are uniformly equivalent. 

In the applications (e.g. for the analysis of the spectrum of Witten Laplacian \cite{HeKlNi04_01} or Kramers-Fokker-Planck operator \cite{HeHiSj11_01}), the supersymmetric structure is used to 
make link between the spectrum of $P$ and the spectrum of the associated operator on $1$-forms. 
As we have seen before, this operator is well-defined if the matrix $G(x;h)$ is invertible. 
It is then important to find condition which insure that $G$ is invertible. We will come back to this issue at the end of section 2.

In practice, it is often useful to expand quantities in powers of the semiclassical parameter $h$. Given a function $f\in S(a)$,
we shall say that it has a classical expansion if there exists a sequence $(f_k)_{k\in\N}$ in $S(a)$ such that for all $K\in\N$,
$$f-\sum_{k=0}^Kh^kf_k\in S(h^{K+1}a).$$
We shall denote by $S_{cl}(a)$ the set of semiclassical functions having a classical expansion and $S_{m,cl}=S_{cl}(\rho_m)$.
Let us now recall one of the results proved in \cite{HeHiSj13}.
\begin{theorem}\label{th:HeHiSj}
Let $P=P(x,hD_x;h)$ be a second order scalar real semiclassical differential operator on $X$. 
Assume there exists $\varphi,\psi\in \ccc^\infty(X,\R)$ such that $P(e^{-\varphi/h})=P^*(e^{-\psi/h})=0$, 
where $P^*$ denotes the formal adjoint of $P$. Assume also that 
the $\delta$ complex is exact in degree $1$ for smooth sections. Then $P$ has a supersymmetric structure.
\end{theorem}
Notice that the above theorem holds true in a very general context. 
For instance, if $X=\R^n$, any scalar second order differential operator such that 
$P(e^{-\varphi/h})=P^*(e^{-\psi/h})=0$ admits a supersymmetric structure. 
Nevertheless, no control on the linear map $G(x;h)$ giving the supersymmetric structure is proved. 
This was noticed in the remark after Definition 1.1 of \cite{HeHiSj11_01} and finding
condition that insure a control over $G(x;h)$ was raised as an open question. Moreover, the author emphasize the fact that the procedure of factorization of $P$ runs in two 
separate ways. The factorization of the symmetric part of $P$ is immediate and provides explicit bound whereas the antisymmetric part is obtained as a solution of a
$\delta$ problem with exponential weights. This difficulty concentrated on the antisymmetric part has to be linked with the general factorization result obtained
for pseudodifferential operators in \cite{BoHeMi15}. 

Let us now recall briefly this last result. We state the result in a slightly different class of symbol that one the used in 
\cite{BoHeMi15} in order that it contains the second order differential operators. It is not difficult to see that the proof could be adapted to the present context.
Let $p(x,\xi)$ be a symbol in the class $S(\<\xi\>^2)$, and let 
$P=\Op_h^w(p)$ denote its Weyl quantization (we refer to \cite{DiSj99_01} for basics of pseudodifferential calculus). Assume that $\varphi$ is a smooth function that behaves as $|x|$ at infinity such that 
$P(e^{-\varphi/h})=0$. The fundamental assumption made in \cite{BoHeMi15} is the following
\be\label{eq:hyp-parity}
\xi\mapsto p(x,\xi)\text{ is even for all }x\in\R^n.
\ee
From this assumption and the equation $P(e^{-\varphi/h})$ we get  $P^*(e^{-\varphi/h})$, and the question of supersymmetry can be investigated.
From Lemma 3.2 and Remark 3.3 in \cite{BoHeMi15}, it follows that 
there exists a matrix-valued pseudodifferential operator $Q_h(x,hD)\in\Psi(1)$ such that 
$P=d_{\varphi,h}^*Q_h(x,hD)d_{\varphi,h}$. 
In other word, $P$ admits a temperate supersymmetric structure with the matrix $G(x;h)$ replaced by the pseudo $Q_h(x,hD)$.

In order to discuss the preceding results and state our first theorem, we need to write the operator $P$ in a specific form. 
It is not hard to verify that any second order scalar real semiclassical differential operator on $X$ can be written in a unique way under the form
\be\label{eq:P-glob}
P(x,hD_x,h)=h\delta \circ A(x;h)\circ hd+U(x;h)\circ hd+ v(x;h)
\ee
where $A$, $U$ and $v$ have the following properties: 
 \begin{itemize}
 \item $A(x;h):T_x^*X\rightarrow T_xX$ and $U(x;h):T_x^*X\rightarrow \R$ are linear and  $v(x;h)\in\R$ 
 \item identifying $T^{**}_xX$ and $T_xX$, $A(x;h)$ is symmetric.
 \item $A$, $U$ and $v$ belong to $S_m$ for some $m\in\R$.
 \end{itemize}
 Observe that $U(x;h)\in T^{**}_xX$ for any $x\in X$. Using again the canonical  identification $T^{**}_xX\simeq T_xX$, it can be seen as an element of $T_xX$.
 In local coordinates, \eqref{eq:P-glob} reads
\be\label{eq:P-loc-coord}
P=-\sum_{i,j=1}^nh\partial_{x_i}\circ a_{i,j}(x;h)\circ h\partial_{x_j}+\sum_{k=1}^n u_k(x;h)\circ h\partial_{x_k}+v(x;h)
\ee
for some real symmetric matrix $A=(a_{ij}(x;h))$, some vector $U=(u_k(x;h))$ and $v(x;h)\in\R$.
So that we can compare the results in \cite{BoHeMi15} and \cite{HeHiSj13}, we shall rewrite the assumptions  in a pseudodifferential way. 
Assume that we work on $X=\R^n$ and that $P$ is given by \eqref{eq:P-glob}, then one has $P=\Op_h^w(p)$, $p=p_{even}+p_{odd}$ with
\be
p_{even}(x,\xi)=\xi^tA(x)\xi+v(x)+\frac h {2}div(U)+\frac{h^2}4\sum_{i,j}\partial_i\partial_ja_{ij}(x)\text{ and }
p_{odd}(x,\xi)=-iU(x)\xi
\ee
where we dropped the dependence of the functions with respect to $h$ in order to lighten the notations.
Suppose that $U(x)=0$ and $P(e^{-\varphi/h})=0$, then we can apply the results of \cite{BoHeMi15} 
and the operator $P$ admits a temperate supersymmetric structure. 
This structure obtained from \cite{BoHeMi15} is a priori a pseudodifferential operator $Q$. 
Nevertheless, a careful look at the proof shows that the operator $Q$ is in fact a multiplication by a temperate matrix.
In the case where $\varphi=\psi$ this gives an improvement of  the conclusion of 
Theorem \ref{th:HeHiSj}. In the following we try to find sharp assumption to make on the antisymmetric part of $P$ in order to prove temperate supersymmetry.

Let us recall the general framework. We consider a second order scalar semiclassical differential operator written under the 
form \eqref{eq:P-glob} and we assume that there 
exists $\varphi,\psi$ such that $P(e^{-\varphi/h})=P^*(e^{-\psi/h})=0$. We wonder if $P$ admits a temperate supersymmetric structure. 
A simple computation shows that   $P(e^{-\varphi/h})=0$ if and only if the following eikonal equation holds true:
\be\label{eq:eik1}
d\varphi\righthalfcup (A(x)d\varphi)+U(x)d\varphi-v(x)+h\delta(A(x)d\varphi)=0.
\ee
On the other hand, since $A$ is symmetric, we have  $P^*=h\delta\circ A(x)\circ hd-U\circ hd+h\delta(U)+v(x)$, 
and we see that $P^*(e^{-\psi/h})=0$ is equivalent to a second eikonal equation:
\be\label{eq:eik2}
d\psi\righthalfcup(A(x)d\psi)-U(x)d\psi-v(x)-h\delta(U)+h\delta(A(x)d\psi)=0
\ee
where $U(x):T_x^*X\rightarrow\R$ is sometimes seen as an element of $T_xX$.

For any  $\phi\in\ccc^\infty(X)$ and $N\in\N$, let 
$E^N_\phi\subset \ccc^\infty(X,\Lambda^2TX)$ denote the subspace of $\ccc^\infty(X,\Lambda^2TX)$ given by 
$$E^N_\phi=\{\alpha=\sum_{finite} (\alpha_j\circ\phi)\, \theta_j,\;\alpha_j\in S(\<t\>^N), \, \theta_j\in ker(\delta)\}$$
and define also the corresponding classical set by
$$E^N_{\phi,cl}=\{\alpha\in E_\phi,\;\alpha_j\in S_{cl}(\<t\>^N)\}.$$
We are now in position to state our first result.
\begin{theorem}\label{th:susy-gen} Let $P$ be as in \eqref{eq:P-glob} with coefficients $A,U,v$ belonging to $S_{m_1}$ for some 
$m_1\in\R$.
Assume that there exists $\varphi,\psi\in S_{m_2}$ for some $m_2\in\R$ such that \eqref{eq:eik1} 
holds true and assume that
\be\label{eq:structop}
U+d(\varphi-\psi)^\righthalfcup\circ A\in\delta(E^N_{\varphi+\psi}).
\ee
Then $P$ admits a  temperate supersymmetric structure given by some $G(x;h):T_x^*X\rightarrow T_xX$.

If additionally, $A$ and $v$ are classical functions and $U+d(\varphi-\psi)^\righthalfcup\circ A\in\delta(E^N_{\varphi+\psi,cl})$, then 
the linear map $G(x;h)$ has a classical expansion.
\end{theorem}
Observe that the conclusion of the above theorem implies that the second  second eikonal equation \eqref{eq:eik2} holds true. 
This could look surprising since we  did not specifically require \eqref{eq:eik2} in our assumptions. In fact, one can easily  prove that 
\eqref{eq:eik1} and \eqref{eq:structop} imply \eqref{eq:eik2}. 
It is natural to wonder if assumption \eqref{eq:structop} is necessary in order to have temperate supersymmetric structure. 
In last section of this note we give partial answer to this question.
 This work was supported by the European Research Council, 
ERC-2012-ADG, project number 320845:  Semi Classical Analysis of Partial Differential
Equations.

\section{Proof of Theorem \ref{th:susy-gen}}
For any antisymmetric $H(x):T^*_xX\rightarrow T_xX$, define $\delta(H):T_xX\rightarrow\R$  by $\delta\circ H\circ d=\delta(H)\circ d$ 
(since $H$ is antisymmetric, this operator is indeed an homogeneous first order differential operator).
For any $G=G(x;h):T_x^*X\rightarrow T_xX$, one has on the $1$-forms
\begin{equation*}
\begin{split}
d_{\psi,h}^{G,*}d_{\varphi,h}&=h\delta\circ G^t\circ hd+d\psi^\righthalfcup\circ G^t\circ hd+h\delta\circ G^t(d\varphi)+
d\psi^\righthalfcup\circ G^t\circ d\varphi^\wedge\\
&=h\delta\circ G^t\circ hd+d\psi^\righthalfcup\circ G^t\circ hd-d\varphi^\righthalfcup\circ G\circ hd+h\delta(G^td\varphi)
+d\psi\righthalfcup G^t(d\varphi)\\
&=h\delta\circ \frac{G^t+G} 2\circ hd+h\frac{\delta(G^t-G)} 2\circ hd+(d\psi^\righthalfcup\circ G^t-d\varphi^\righthalfcup\circ G)\circ hd\\
&\phantom{++++++++++++++++}+h\delta(G^td\varphi)+d\psi\righthalfcup G^t(d\varphi)
\end{split}
\end{equation*}
Let us introduce the symmetric and antisymmetric part of $G$:
$$G^s=\frac 1 2(G+G^t)\text{ and }G^a=\frac 12 (G-G^t)$$
then we get 
\be\label{eq:devt-susy}
\begin{split}
d_{\psi,h}^{G,*}d_{\varphi,h}=h\delta\circ G^s\circ hd&-h\delta(G^a)\circ hd+(d(\psi-\varphi)^\righthalfcup\circ G^s-d(\varphi+\psi)^\righthalfcup\circ G^a)\circ hd\\
&+h\delta(G^sd\varphi)-h\delta(G^ad\varphi)+d\psi\righthalfcup(G^sd\varphi)-d\psi\righthalfcup(G^ad\varphi)
\end{split}
\ee
Identifying \eqref{eq:devt-susy} and \eqref{eq:P-glob}, we see that $P=d_{\psi,h}^{G,*}d_{\varphi,h}$ if and only if 
\be\label{eq:ident-susy1}
\left\{
\begin{array}{c}
 G^s(x)=A(x;h)\phantom{+++++++++++++}\\
 U(x;h)+d(\varphi-\psi)^\righthalfcup \circ G^s=-h\delta(G^a)-d(\varphi+\psi)\righthalfcup G^a\phantom{+}\\
 v(x;h)=h\delta(G^sd\varphi)-h\delta(G^ad\varphi)+d\psi\righthalfcup(G^sd\varphi)-d\psi\righthalfcup(G^ad\varphi)
\end{array}
\right.
\ee
Looking for $G$ under the form $G=A+B$ with $B$ antisymmetric, \eqref{eq:ident-susy1} becomes
\be\label{eq:ident-susy2}
\left\{
\begin{array}{c}
 U(x;h)+d(\varphi-\psi)^\righthalfcup\circ A=-h\delta(B)-d(\varphi+\psi)^\righthalfcup\circ B\phantom{+}\\
 v(x;h)=h\delta(Ad\varphi)-h\delta(Bd\varphi)+d\psi\righthalfcup(Ad\varphi)-d\psi\righthalfcup(Bd\varphi)
\end{array}
\right.
\ee
Suppose now that the first equation of the above system is solved. 
Then, 
$$U(d\varphi)=d\psi\righthalfcup(Ad\varphi)-d\varphi\righthalfcup(Ad\varphi)-h\delta(Bd\varphi)+d\psi\righthalfcup Bd\varphi$$
and using \eqref{eq:eik1}, we get easily the second one.
Hence, we are reduced to find a map $B\in \ccc^\infty(X,\Lambda^2 TX)$ which is temperate and solves
\be\label{eq:ident-susy3}
U+d(\varphi-\psi)^\righthalfcup\circ A=-h\delta(B)-d\phi\righthalfcup B
\ee
where $\phi=\varphi+\psi$. Thanks to assumption \eqref{eq:structop}, there exists $\theta_1,\ldots,\theta_K\in \ccc^\infty(X,\Lambda^2TX)$ and 
$\alpha_1,\ldots\alpha_K\in\ccc^\infty(\R,\R)$ such that $\delta \theta_k=0$ for all $k$ and 
$$U+d(\varphi-\psi)^\righthalfcup\circ A=\delta\theta$$
with $\theta=\sum_{k=1}^K(\alpha_k\circ\phi)\,\theta_k$. Hence \eqref{eq:ident-susy3} is equivalent to
\be\label{eq:ident-susy4}
\delta\theta=-h\delta(B)-d\phi\righthalfcup B.
\ee
On the other hand, for any $k$, we have 
$$\delta ((\alpha_k \circ\phi)\theta_k)=(\alpha_k \circ\phi)\delta\theta_k-d(\alpha_k \circ\phi)\righthalfcup \theta_k=-(\alpha'_k\circ\phi)d\phi\righthalfcup\theta_k$$ and hence 
\eqref{eq:ident-susy4} is equivalent to 
\be\label{eq:ident-susy5}
h\delta(B)+d\phi\righthalfcup B=\sum_{k=1}^K(\alpha'_k\circ\phi )d\phi\righthalfcup\theta_k.
\ee
Since this is a linear equation, it suffices find some temperate $B_k$ such that 
$$
h\delta(B_k)+d\phi\righthalfcup B_k=(\alpha'_k\circ\phi )d\phi\righthalfcup\theta_k.
$$
for all $k$. In order to lighten the notation, we will drop the index $k$ in the following lines. Setting $\tilde B=e^{-\phi/h} B$, the above equation is equivalent to
\be\label{eq:ident-susy6}
h\delta(\tilde B)=e^{-\phi/h}(\alpha'\circ\phi )d\phi\righthalfcup\theta.
\ee 
Our aim is to find a solution $\tilde B$ of this equation such that $B=e^{\phi/h}\tilde B$ is temperate.
For this purpose, simply observe that 
$$e^{-\phi/h}(\alpha'\circ\phi )d\phi=d(\beta\circ\phi)$$
with $\beta\in \ccc^\infty(\R,\R)$ defined by
\be\label{eq:def-beta}
\beta(t)=-\int_t^{m_\infty}\alpha'(s)e^{-s/h}ds
\ee
with $m_\infty=+\infty$ when $X=\R^n$ and $m_\infty=\sup \phi$ when $X$ is a compact manifold.
Hence \eqref{eq:ident-susy6} becomes
$$
h\delta(\tilde B)=d(\beta\circ\phi)\righthalfcup\theta=-\delta((\beta\circ\phi)\theta).
$$
A solution is trivially given by $\tilde B=-\frac 1 h(\beta\circ\phi)\theta$, that is 
\be
\begin{split}
B(x)&=\big(\frac 1 he^{\phi/h}\int_{\phi(x)}^{m_\infty}\alpha'(s)e^{-s/h}ds\big) \theta(x)\\
&=\frac 1 h\big(\int_0^{m_\infty-\phi(x)}\alpha'(s+\phi(x))e^{-s/h}ds\big) \theta(x)
\end{split}
\ee
It remains to check that $B$ is temperate. For this purpose it suffices to observe that $m_\infty-\phi(x)\geq 0$ and hence,
we have necessarily $e^{-s/h}\leq1$ in the above integral.
In the case where $X$ is compact, this shows immediately that $\partial^\nu B=\ooo(\rho_{Nm_2})$ for all $\nu\in\N^n$. 
In the case where $X=\R^n$, using the fact that $\alpha$ has at most polynomial growth and
performing integration by parts we obtain similarly $\partial^\nu B=\ooo(\rho_{Nm_2})$. 

Suppose now  that $A,v$ have classical expansion and that $U+d(\varphi-\psi)^\righthalfcup\circ A\in\delta(E^N_{\varphi+\psi,cl})$.
In order to show that $G$ has classical expansion, it suffices to do so for $B(x)$ above.
A simple change of variable shows that 
$$B(x)=-\big(\int_0^{(m_\infty-\phi(x))/h}\alpha'(hs+\phi(x))e^{-s}ds\big) \theta(x)$$
In the case where $X=\R^n$, $m_\infty=+\infty$ and a simple Taylor expansion in the above integral gives the result.
Suppose now that $X$ is a compact manifold, Taylor expansion of $\alpha'(hs+\phi(x))$ reduce the proof to 
expand terms of the form
$$\int_0^{(m_\infty-\phi(x))/h}(hs)^ke^{-s}ds$$ which is easily obtained by integration by parts.
This concludes the proof.
\ep
\begin{remark}
 From the preceding proof, we see that if $\alpha$ satisfies additional properties, then $B$ can be computed explicitly. 
 For instance, if $\alpha$ is polynomial, integration by parts leads to explicit formula for $B$.
\end{remark}

In conclusion of this section we shall discuss the invertibility of $G(x;h)$. In order to simplify the discussion, we suppose that $X=\R^n$.
Assume that there exists $\theta\in E_{\varphi+\psi}$ such that  $U+d(\varphi-\psi)^\righthalfcup\circ A=\delta\theta$. 
From the proof above, one has 
$G=A+B$ with $A$ real symmetric defined by \eqref{eq:P-glob} and $B$ real antisymmetric.
Assume that $A$ is uniformly positive definite, that is
$$\exists C>0,\;\forall x\in \R^n,\,\forall\xi\in T_x^*X,\;\<A(x;h)\xi,\xi\>\geq C\vert \xi\vert^2.$$
One check easily that $G$ enjoys the same estimate and hence is invertible with $G^{-1}$  bounded by $C^{-1}$.

Let us now  consider the case where  $A(x;h)$ is only positive (not necessary definite). 
Assume that there exists some orthogonal subspaces  $E,F$ independent of $(x,h)$ such that 
$\R^n=E\oplus F$ with $E=\ker A(x;h)$ for all $x\in\R^n$. In some specific situations we can insure that $G$ is invertible.
For instance, if $ker B=F$ and $B_{|E}$ and $A_{|F}$ are invertible with inverse uniformly bounded, then $G$ has a uniformly bounded inverse.
Another interesting situation is a generalization of Kramers-Fokker-Planck operator. Recall that the Kramers-Fokker-Planck operator is defined
on $\R^{2n}$ by
\be\label{eq:KFP}
K(h)=y\cdot h\nabla_x-\nabla_xV(x)\cdot h\nabla_y-h^2\Delta_y+y^2-hn
\ee
where $(x,y)\in\R^n\times\R^n$ denotes the space variable and $V$ is a smooth function. This operator admits a supersymmetric structure
$K(h)=d_{\varphi,h}^{G,*}\circ d_{\varphi,h}$ with $\varphi(x,y)=\frac 12 y^2+V(x)$ and 
$$G(x,h)=\frac 12\left(
\begin{array}{cc}
 0&\Id\\
 -\Id& 2 \Id
 \end{array}
\right).
$$
In particular, $G$ is invertible and $G^{-1}=\ooo(1)$. In \cite{HeHiSj11_01}, this permitted the authors to define the whole associated Witten complex and 
to obtain precise information on the spectrum of $K(h)$.
This situation can be easily generalized. Let us go back to the above situation where  $E=\ker A(x,h)$ is independent of $x$ and 
suppose  that $B(E)\subset F$. Denote $\Pi:X\rightarrow E$ the orthogonal projection and $\hat\Pi=1-\Pi$. Then, $\Pi B\Pi=0$ and
$$G=\hat \Pi A\hat \Pi+\Pi B\hat\Pi+\hat\Pi B\Pi+\hat\Pi B\hat \Pi.$$
Therefore, the equation $G\xi=0$ is equivalent to 
$$
\left\{
\begin{array}{cc}
 \hat \Pi A\hat \Pi\xi+\hat\Pi B\Pi\xi +\hat\Pi B\hat \Pi\xi=0\\
\Pi B\hat\Pi\xi=0.
 \end{array}
\right.
$$
Taking the scalar product with $\xi$, we get $\<\hat \Pi B\Pi\xi,\xi\>=\<\Pi\xi,B^t\hat \Pi\xi\>=-\<\xi,\Pi B\hat\Pi\xi\>=0$ and hence
$$\<A\hat \Pi\xi,\hat\Pi\xi\>+\<B\hat\Pi\xi,\hat\Pi\xi\>=0.$$
Since $A$ is definite positive on $F$ and $B$ is antisymmetric it follows immediately that $\hat\Pi\xi=0$. Hence 
$G$ is injective if and only if $B_{|E}:E\rightarrow F$ is injective. In particular it is necessary that $\dim F\geq \dim E$.
In order to estimate $G^{-1}$ one starts  from $G\xi=\eta$. Working as above we show easily that $\|\hat\Pi\xi\|\leq 2C_0^{-1}\|\eta\|$ where 
$C_0$ is the smallest eigenvalue of $A$ on $F$. One has also
$$\hat\Pi B\Pi\xi=\hat\Pi\eta-\hat \Pi A\hat \Pi\xi-\hat\Pi B\hat \Pi\xi.$$
Hence, assuming additionally that $B:E\rightarrow Im(B)$ has a bounded inverse with respect to $h$, 
 $G^{-1}$ is automatically bounded with respect to $h$.

\section{About assumption \eqref{eq:structop}}

The aim of this section is to discuss the necessity of assumption \eqref{eq:structop} in order to get temperate supersymmetry. 
Throughout this section, we assume  that $\varphi=\psi$. Then, the question of the existence of a supersymmetric structure can be reduced 
to the same question for operators of order $1$.
Indeed, we can write $P=P_1+P_2$ with $P_1=\frac 1 2(P+P^*)$. Observe that $P_1$ is formally selfadjoint and $P_2$ is formally antiadjoint. 
Moreover, since $\varphi=\psi$, we have
$P_1(e^{-\varphi/h})=0$, $P_2(e^{-\varphi/h})=0$ and we claim that $P_1$ admits automatically a temperate supersymmetric structure. Indeed, we have
$$P_1=h\delta\circ A\circ hd+\frac h 2\delta(U)+v$$
and we know from the proof of Theorem \ref{th:susy-gen}, that $P_1$ admits a temperate supersymmetric structure if and only if 
there exists $B$ antisymmetric, such that 
$$
d(\varphi-\psi)^\righthalfcup\circ A=-h\delta(B)-d(\varphi+\psi)^\righthalfcup\circ B$$
Since $\varphi=\psi$, $B=0$ solves  this equation.
We are then reduced to investigate the condition that insure that 
\be\label{eq:reduc_op_1}
P_2=U\circ hd-\frac h 2\delta(U)
\ee
admits  a temperate supersymmetric structure. 

\subsection{The two dimensional case}
In this section we assume that $X=\R^2$ the euclidean plane and we consider operators $P$ of the form \eqref{eq:reduc_op_1} in the case where $\delta(U)=0$. 
We denote by $\mathscr{C}$ the set of critical points of $\varphi$ and by $w$ the canonical $2$-form. Using the euclidean structure we can write $P=U(x)\cdot h\nabla$. 
The following lemma shows that away from critical points $U$ has necessarily the form \eqref{eq:structop}.
\begin{lemma}\label{lem:reduc_champ-2D}
For any $x_0\in\R^2\setminus\mathscr{C}$, there exists a neighborhood $V$ of $x_0$ and a smooth function $f_{x_0}:\R\rightarrow\R$ such that 
$U=\delta ((f_{x_0}\circ\varphi) \omega)$.
\end{lemma}
\bp
Under the above assumptions, the eikonal equation reads
$$U\cdot\nabla\varphi=0$$
and since $\delta(U)=0$, there exists $\alpha\in \ccc^\infty(\R^d)$ such that $U=\delta(\alpha dx_1\wedge dx_2)=\partial_2\alpha dx_1-\partial_1\alpha dx_2$ where we identify $U$ with a $1$-form.
Going back to the eikonal equation we obtain $\partial_2\alpha \partial_1\varphi-\partial_1\alpha\partial_2\varphi=0$ which can be interpreted as $\det(\nabla\alpha,\nabla\varphi)=0$.
From this equation we deduce that near any point $x_0\in\R^2$ there exists a smooth function $f_{x_0}$ such that $\alpha=f_{x_0}\circ\varphi$. 
To see this, recall  that $x_0$ in non-critical, hence there exists $V$  neighborhood of $x_0$ such that $V\cap\mathscr{C}=\emptyset$. Shrinking $V$ if necessary and changing coordinates, we can assume that there exists $\nu\in \R^2$ such that 
\be\label{eq:redress-phi1}
\varphi(z)=\varphi(x_0)+\<\nu,z-x_0\>\text{ for all }z\in V.
\ee
Suppose that $x,y\in V$ are such that $\varphi(x)=\varphi(y)$. 
Then, we deduce  from \eqref{eq:redress-phi1} that there exist a smooth path $\gamma:[0,1]\rightarrow V$ such that $\gamma(0)=x$, $\gamma(1)=y$ and $\varphi\circ\gamma$ is constant
(take $\gamma(t)=x+t(y-x)$). Since
$$\alpha(x)-\alpha(y)=\int_0^1\frac d{dt}\alpha\circ\gamma(t)dt=\int_0^1\nabla\alpha\circ\gamma\cdot\dot\gamma(t)dt=0$$
for that $\dot\gamma $ is orthogonal to $\nabla\varphi$, this shows that $\alpha$ depends only on $\varphi$ on $V$. Hence there exists a function $f_{x_0}$ such that $\alpha=f_{x_0}\circ\varphi$.
Moreover, using the fact that $\nabla\varphi$ doesn't vanish on $V$ one can easily show that $f_{x_0}$ is smooth.
\ep

Without additional assumption on $\varphi$ it seems difficult to globalize the above result. However, in the case where $\varphi$ is a Morse function 
one can get further information on the structure of $U$. In the following, we assume that $\varphi$ is a Morse function  such that 
 $\CC$ is finite.
For $k=0,1,2$, let $\CC^{(k)}$ denote the set of critical points of index $k$. 
The following is an improvement of Lemma \ref{lem:reduc_champ-2D}.
\begin{lemma}\label{lem:reduc_champ-2D-Morse}
Assume that $\varphi$ is a Morse function, then  for any $x_0\in \R^2\setminus\CC^{(1)}$ there exists a neighborhood $V$ of $x_0$ and a smooth function $f_{x_0}:\R\rightarrow\R$ such that 
$U=\delta ((f_{x_0}\circ\varphi) \omega)$.
\end{lemma}
\bp
It suffices to check the conclusion in the case where $x_0$ in either a minimum or a maximum of $\varphi$. Assume that $x_0$ is a minimum of $\varphi$ (the maximum case can be treated in the same way). As in the previous lemma, we first choose a small neighborhood $V$ of $x_0$ and new coordinates such that 
\be\label{eq:redress-phi2}
\varphi(z)=\varphi(x_0)+\vert z-x_0\vert^2\text{ for all }z\in V
\ee
Without loss of generality, we can also assume that $x_0=0$ and $\varphi(0)=0$.
Let $x,y\in V$ be such that $\varphi(x)=\varphi(y)$, that is $\vert x\vert=\vert y\vert$ and denote by $\alpha$ the angle between $x$ and $y$.
Let  $\gamma:[0,1]\rightarrow V$ be the path defined by $\gamma(t)=r_{t\alpha}(x)$, where $r_\theta$ denotes the rotation of angle $\theta$. Then 
$\gamma(0)=x$, $\gamma(1)=y$ and $\varphi\circ\gamma$ is constant. The same argument as in Lemma \ref{lem:reduc_champ-2D} shows that $\alpha(x)=\alpha(y)$.
Hence $\alpha$ depends only on $\varphi$ on $V$ and there exists a function $f_{x_0}$ such that $\alpha=f_{x_0}\circ\varphi$. It is clear that $f_{x_0}$ is smooth away from $x_0=0$ since the gradient 
of $\varphi$ doesn't vanish. In order to show that $f_{x_0}$ is smooth in $x_0$ we write 
$f_{x_0}(t)=\alpha(\sqrt t,0)$. Let us write the Taylor expansion of $\alpha$ near the origin
$$\alpha(x_1,x_2)\simeq\sum_{j,k}\alpha_{j,k}x_1^jx_2^k.$$
The equation $\det(\nabla\alpha,\nabla\varphi)=0$ yields $x_2\partial_1\alpha=x_1\partial_2\alpha$ and it follows that for all $j,k$, $j\alpha_{j,k-2}=k\alpha_{j-2,k}$ with the convention 
$\alpha_{p,q}=0$ for $p<0$ or $q<0$. Using this relation, we get immediately $\alpha_{j,k}=0$ for any $j,k$ such that $j$ or $k$ is odd. In particular $\alpha_{j,0}=0$ for any odd $j$, which shows that 
$f_{x_0}$ is smooth at the origin.
\ep

Let $\Sigma=\varphi(\CC^{(1)})\subset\R$ denote the saddle values of $\varphi$. 
Then $\R^2\setminus \varphi^{-1}(\Sigma)$ has a finite number of connected components $\Omega_1,\ldots,\Omega_J$ and only one (say $\Omega_J$) is unbounded.
\begin{lemma}\label{lem:reduc_champ-2D-Morse-global}
Assume that  $\varphi$ is a Morse function, then there exists some functions $f_1,\ldots,f_J\in\ccc^\infty(\R)$  such that 
$$ U=\delta ((f_j\circ\varphi) \omega) \text{ on }\Omega_j$$
\end{lemma}
\bp Let $j\in \{1,\ldots,J\}$ be fixed and 
let $x,y\in \Omega_j$ such that $\varphi(x)=\varphi(y)=:\sigma$. We shall prove that the set $F_\sigma:=\varphi^{-1}(\sigma)$ is arcwise connected. 
We first assume that $\Omega_j$ is bounded. Hence there exists a covering of $\overline\Omega_j$ by a finite collection of convex  open set $(\omega_k)_{k=1,\ldots, K}$ such such that on each 
$\omega_k$, there exists some change of coordinates $\theta_k:\ooo_k\rightarrow\omega_k$ such that $\varphi_k=\varphi\circ\theta_k$ takes one of the following forms:
$$\varphi_k(z)=\<z,\nu\>,\,\nu\in\R^2\setminus 0\text{ or }\varphi_k(z)=\vert z\vert^2\text{ or }\varphi_k(z)=-\vert z\vert^2$$
for any $z\in\ooo_k$ neighborhood of $0\in\R^2$.
Let $M=2\pi K\sup_{k=1,\ldots,K}(\Vert D\theta_k\Vert_\infty diam(\ooo_k))$ and 
\be\label{eq:def-Gamma}
\Gamma_j=\{M-\text{Lipschitz path }\gamma\text{ contained in }\Omega_j\text{ and joining } x \text{ to } y\}.
\ee
 Since $\Omega_j$ is arcwise connected, $\Gamma_j$ is nonempty. Indeed
there exists a smooth path $\gamma_0:[0,1]\rightarrow\Omega_j$ joining $x$ to $y$,  and up to reparametrization we can also assume that $\vert\gamma_0'\vert$ is constant. Moreover, using the specific form 
of $\varphi$ on each $\omega_k$ we can modify $\gamma_0$ into a piecewise $\ccc^1$ path  so that $I_k:=\{t\in[0,1],\,\gamma_0(t)\in\omega_k\}$
is an interval for all 
$k=1,\ldots,K$. It follows easily that 
$$|\gamma_0'(t)|\leq 2\pi\sum_{k,\;I_k\neq\emptyset}diam(\ooo_k)\Vert D\theta_k\Vert_\infty\leq M$$
excepted for a finite number of values of $t$. Therefore  $\gamma_0\in\Gamma_j$. 

Introduce next the set $\mmm=\{\sup_{[0,1]}\varphi\circ\gamma,\;\gamma\in\Gamma_j\}\subset[\sigma,+\infty[$ since $\varphi(x)=\sigma$, 
and let $m=\inf\mmm\geq \sigma$. We claim that $m=\sigma$.
Indeed, if follows from Ascoli theorem that $\Gamma_j$  is relatively compact in $\ccc([0,1],\overline\Omega_j)$. 
Hence there exists a path $\gamma_1$ contained in $\Omega_j$ and joining $x$ and $y$ such that $m=\sup\varphi\circ\gamma_1=\varphi\circ\gamma_1(t_1)$. 
Suppose by contradiction that $m>\sigma$ and let $x_1=\gamma_1(t_1)$. 
By definition of $\Omega_j$, $x_1$ can not be a saddle point of $\varphi$ and since $m=\sup\varphi\circ\gamma_1$ it  is neither  a minimum. 
Hence $x_1$ is either a local maximum either a noncritical point of $\varphi$. 
In both case, it is easy to modify locally the path $\gamma_1$ in order to decrease $m$. This gives a contradiction.

Hence $m=\sigma$ and there exists a continuous path $\tilde\gamma_1\subset \Omega_j$ 
joining $x$ and $y$ and such that $\sup_{[0,1]}\varphi\circ\tilde\gamma=\sigma$. Moreover, by construction $\tilde\gamma_1$ is $M$-Lipschitz.
Therefore, the set 
\be\label{eq:def-tildeGamma}
\tilde\Gamma_j=\{M-\text{Lipschitz path }\gamma\text{ contained in }\Omega_j\cap\{\varphi\leq\sigma\}\text{ and joining } x \text{ to } y \}.
\ee 
is nonempty. Let $\lll=\{\inf_{[0,1]}\varphi\circ \tilde\gamma,\,\gamma\in \tilde\Gamma_j\}$ and 
$\ell=\sup\lll$. As before, there exists  a Lipschitz path $\gamma_2$ such that $\ell=\sup\varphi\circ\gamma_2$ and we can show easily that $\varphi\circ\gamma_2$ is constant equal to $\sigma$.
Using this path $\gamma_2$ and the fact that $\varphi^{-1}(\sigma)$ is locally connected, we construct a path $\gamma_3\subset \varphi^{-1}(\sigma)$ from $x$ to $y$ which is piece wised $C^1$.

Using this  path $\gamma_3$ and repeating the argument of the proof of Lemma \ref{lem:reduc_champ-2D-Morse}, it follows easily that $\alpha(x)=\alpha(y)$ and hence $\alpha$ depends only on $\varphi$. 
This permits to construct a function $f_j$ such that $\alpha=f_j\circ \varphi$ on $\Omega_j$. The smoothness of $f_j$ is a local property and then follows from Lemma \ref{lem:reduc_champ-2D-Morse}.

Let us now prove the result for the unbounded component $\Omega_J$.
Let $\sigma\in\R$ be fixed and let $x,y\in \Omega_J\cap \varphi^{-1}(\sigma)$. By definition, there exists a path $\gamma$ contained in $\Omega_J$ joining 
$x$ and $y$. Let $R>0$ be such that $\gamma\subset B(0,R)$. Since $\Omega_J\cap B(0,R)$ is relatively compact, one can follow the same strategy as for the 
bounded component with $\Omega_j$ replaced by $\Omega_J\cap B(0,R)$.
\ep

As a consequence of the above lemma we get the following

\begin{theorem} Let $P(h)=U\circ hd$ with $\delta(U)=0$. Assume that $\varphi$ is a Morse function with a finite number of critical points and such that $U\cdot\nabla\varphi=0$.
Assume additionally that for all $i,j=1,\ldots,J$, $i\neq j$ and all $x\in\Omega_i, y\in\Omega_j$ such that $\varphi(x)=\varphi(y)$ there exists a smooth path $\gamma$ from $x$ to $y$ such that 
\be\label{eq:lien_cc}
\int_\gamma \star U=0
\ee
 where $\star$ denotes the Hodge star operator. Then $P$ satisfies  \eqref{eq:structop} and hence admits a temperate supersymmetric structure.
\end{theorem}
\bp
Let $I$ denote the image of $\varphi$ which is a (bounded or unbounded) interval.
From \eqref{eq:lien_cc}, one knows that for all $i\neq j$ and all $x\in\Omega_i$, $y\in\Omega_j$ such that $\varphi(x)=\varphi(y)$ one has 
$f_i\circ\varphi(x)=f_j\circ\varphi(y)$. Hence, the function $f:I\rightarrow\R$ given by $f\circ\varphi(x)=\alpha(x)$ is well defined.
One has to show that $f$ is smooth and the only point which has not been already examined is the smoothness near saddle points. 
Let $s_0$ be a saddle point of $\varphi$. Without loss of generality, we assume $s_0=0$ and $\varphi(x)=x_1^2-x_2^2$ near the origin.
As in the proof of Lemma \ref{lem:reduc_champ-2D-Morse}, we write $U=\delta (\alpha\omega)$ with 
$\alpha\simeq\sum_{j,k}\alpha_{j,k}x_1^jx_2^k$, and it follows from the equation $U(d\varphi)=0$ that 
\be\label{eq:aux1}
j\alpha_{j,k-2}=-k\alpha_{j-2,k}
\ee
for all $j,k$ (with the convention $\alpha_{j,k}=0$ for negative $j$ or $k$). As before, we get $\alpha_{j,k}=0$ for $j$ or $k$ odd and one has
$$
f(t)=
\left\{
\begin{array}{cc}
 \alpha(\sqrt t,0)&\text{ if } t>0\\
 \alpha(0,\sqrt{-t})&\text{ if } t<0
\end{array}
\right.
$$
From the Taylor expansion of $\alpha$ we see that $f$ is smooth if and only if $\alpha_{2j,0}=(-1)^j\alpha_{0,2j}$ which is a consequence of \eqref{eq:aux1}.
\ep

In consideration of the above theorem, one could think that operators admitting temperate supersymmetric structure
are more or less of the form \eqref{eq:structop}. 
In fact, when the dimension is greater than $2$, this is not the case.
One way to see this is to notice that supersymmetric structure can be easily tensorized.

Let $X_j$, $j=1,2$ be either an euclidean space either a smooth connected compact manifold. Let $P_j(x_j,hD_{x_j})$ 
denote a second order semiclassical differential operator on $X_j$, and let 
$\varphi_j,\psi_j\in \ccc^\infty(\R)$.
\begin{theorem}\label{th:tensor}
 Assume that the $P_j$ admit a supersymmetric structure $G_j(x_j;h)$ associated to the phase $\varphi_j,\psi_j$. Then 
 the operator $P(x,hD_x)=P_1(x_1,hD_{x_1})+P_2(x_2,hD_{x_2})$ acting on $X=X_1\times X_2$, admits a supersymmetric structure 
 $$P=d_{\psi,h}^{G,*}d_{\varphi,h}$$
 with $\varphi(x_1,x_2)=\varphi_1(x_1)+\varphi_2(x_2)$, $\psi(x_1,x_2)=\psi_1(x_1)+\psi_2(x_2)$ and 
 $G(x;h)=\left(
\begin{array}{cc}
 G_1&0\\
 0&G_2
\end{array}
\right)
$
\end{theorem}
\bp
This is immediate.
\ep

Using this result, we can construct easily some examples where $U$ has not necessarily the form \eqref{eq:structop}, even locally.
For instance, let $X=\R^3$, $\varphi(x)=\psi(x)=x_1^2+x_2^2+x_3^2$ and 
$$P(x,hD_x)=x_1\cos(x_1^2+x_2^2) h\partial_{x_2}-x_2\cos(x_1^2+x_2^2)h\partial_{x_1}=U\circ hd$$ with 
$U(x)=-x_2\cos(x_1^2+x_2^2) dx_1+x_1\cos(x_1^2+x_2^2)dx_2$. Then $P$ admits a supersymmetric structure and one has 
$U=\delta(\frac  12\sin(x_1^2+x_2^2) dx_1\wedge dx_2)$ but $U$ can not be written under the form \eqref{eq:structop}. 
\subsection{Perturbation of supersymmetric structure}

In this section we show that assumption \eqref{eq:structop} is not necessary in general  to get supersymmetry. 

We go back to the general situation where $X$ is either a compact manifold without boundary either $\R^n$. 
We assume that  $\varphi:X\rightarrow\R$ a smooth function such that the set $\vvv=\varphi(\CC)$ of critical values of $\varphi$ is finite.
In the case where $X=\R^n$ we assume additionnaly that $\lim_{|x|\rightarrow\infty}\varphi(x)=\infty$.
For any $\sigma\in\R$ we denote $X_\sigma=\{x\in X,\,\varphi(x)<\sigma\}$ and we consider a fixed connected component $\omega_\sigma$ of $X_\sigma$. For any 
$\epsilon>0$ we denote $\omega_\sigma^\epsilon=\omega_\sigma\cap X_{\sigma-\epsilon}$.
Since $\vvv$ is finite, there exists $\epsilon_0>0$  small enough such that $]\sigma-\epsilon_0,\sigma[\cap\vvv=\emptyset$.
Therefore, the set $\omega_\sigma^\epsilon$ have smooth boundary for all $0<\epsilon<\epsilon_0$ and $\omega_\sigma^\epsilon$ 
is relatively compact in $\omega_\sigma$ 
(in the case $X=\R^n$ this is true  since $\varphi$ goes to infinity at infinity).
Hence we can construct a smooth cut-off function $\chi_\epsilon$ such that $\chi=1$ on $\omega_\sigma^\epsilon$ and $\supp(\chi_\epsilon)\subset \omega_\sigma$.
Let $\alpha\in \ccc^\infty(\R)$ such that $\supp(\alpha)\subset ]-\infty,\sigma-\epsilon[$ and let $\theta $ be a $2$-form such that $\delta\theta=0$.
Consider $$U_\epsilon=\delta((\chi_\epsilon \,\alpha\circ\varphi )\theta).$$
Then $U_\epsilon$ is a smooth $1$-form such that $\delta U_\epsilon=0$. Moreover, by construction 
$\supp(d\chi)\subset \{\sigma-\epsilon<\varphi<\sigma\}$ and hence $\alpha\circ\varphi d\chi=0$. Therefore, we have in fact
$$U_\epsilon=\chi_\epsilon\delta((\alpha\circ\varphi)\theta)=-\chi_\epsilon\alpha'\circ\varphi d\varphi\righthalfcup\theta.$$
and hence $U_\epsilon(d\varphi)=0$.
Then, the operator $P_{\epsilon}=U_{\epsilon}\circ hd$ is formally self-adjoint and we have $P_\epsilon(e^{-\varphi/h})=0$.
\begin{proposition}
Assume that the above assumptions are fulfilled then $P_{\epsilon}$ admits a temperate supersymmetric structure.
\end{proposition}
\bp Set $\phi=2\varphi$.
Let $$B^0(x):=\big(\frac 1 {2h}e^{\phi/h}\int_{\phi(x)}^{2(\sigma-\epsilon)}\alpha'(\frac s 2)e^{-s/h}ds\big) \theta(x)$$
which is temperate since $\alpha$ i supported in $\{s<\sigma-\epsilon\}$.
It follows from the proof of Theorem \ref{th:susy-gen} that 
$$h\delta B^0+d\phi\righthalfcup B^0=(\alpha'\circ\varphi )d\varphi\righthalfcup\theta=-\delta(\alpha\circ\varphi \theta).$$
Set $B=\chi_\epsilon B^0$. Thanks to the support properties of $\chi_\epsilon$ and $\alpha$, one has
$\delta B=\chi\delta B^0$.
Therefore,
$$h\delta B+d\phi\righthalfcup B=-U_\epsilon$$ 
which is exactly the eikonal equation we  have to solve.
Moreover, the same proof as in Theorem \ref{th:susy-gen} with $m_\epsilon$ instead of $m_\infty$ shows that $B$ is temperate.
\ep
\begin{remark}
 Assume $X_\sigma$ has two distinct connected components. 
 Then $U_{\epsilon}$ has a temperate supersymmetric structure and doesn't satisfy
 \eqref{eq:structop}.
\end{remark}

\bibliographystyle{amsplain}
\bibliography{ref_susy}

\end{document}